%% file: main.tex
\documentclass[12pt]{amsart}

\usepackage[left=15mm, right=15mm, top=15mm, bottom=15mm]{geometry}

\usepackage{float}
\usepackage[font=small,labelsep=period, figurename={Рис.}]{caption}

\input{eng}

\input{eng_theorems}
\input{packages}
\input{preamble}

\begin{document}

\input{mains/banach-povm}

\printbibliography

\end{document}

%% file: eng.tex
\usepackage[english]{babel}
\usepackage[T1]{fontenc}
\usepackage[utf8]{inputenc}


%% file: eng_theorems.tex
\theoremstyle{plain}
    \newtheorem{theorem}{Theorem}[section]
    \newtheorem{lemma}[theorem]{Lemma}
    \newtheorem{proposition}[theorem]{Proposition}

\theoremstyle{definition}

    \newtheorem{example}[theorem]{Example}

\usepackage[
    draft = false,
    unicode = true,
    colorlinks = true,
    allcolors = blue,
    hyperfootnotes = true,
    citecolor = red
]{hyperref}

%% file: packages.tex
\usepackage{amsfonts,amsmath,amssymb,amscd,amsthm,url}
\usepackage[inline]{enumitem}
\usepackage{graphicx,epsf,afterpage, wrapfig}
\usepackage{soul}
\usepackage{multicol}
\usepackage{array}
\usepackage{braket}
\usepackage{epigraph}
\usepackage{rotating, floatflt}
\usepackage{xstring} 

\usepackage[
backend=biber,
style=numeric-comp,
sorting=none,
giveninits=true
]{biblatex}
\usepackage{csquotes}

\usepackage{xcolor}

\usepackage{ulem}

\usepackage{tikz}
\usetikzlibrary{positioning, calc, intersections, through, arrows, matrix, chains, math}
\usetikzlibrary{arrows.meta}
\usetikzlibrary{shadows}
\usetikzlibrary{graphs}
\usetikzlibrary{graphs.standard}
\usetikzlibrary{patterns}
\usetikzlibrary{decorations.pathreplacing,calligraphy,backgrounds}

%% file: preamble.tex
\DeclareFieldFormat{usera}{\href{https://arxiv.org/abs/#1}{\texttt{arXiv:#1}}}

\renewbibmacro*{finentry}{\printfield{usera}\newunit\finentry}

\renewbibmacro{in:}{%
  \iffieldundef{usera}
    {\printtext{\bibstring{in}\intitlepunct}}
    {}%
}

\newcommand\norm[1]{\ensuremath{\left\lVert#1\right\rVert}}
\newcommand\abs[1]{\ensuremath{\left\lvert#1\right\rvert}}

\DeclareMathOperator{\Img}{Im}

\DeclareMathOperator{\supp}{supp}

\newcommand{\Acal}{\mathcal{A}}
\newcommand{\Bcal}{\mathcal{B}}

\newcommand{\Fcal}{\mathcal{F}}

\newcommand{\Hcal}{\mathcal{H}}
\newcommand{\Ical}{\mathcal{I}}

\newcommand{\Lcal}{\mathcal{L}}
\newcommand{\Mcal}{\mathcal{M}}

\newcommand{\Rcal}{\mathcal{R}}
\newcommand{\Scal}{\mathcal{S}}

\newcommand{\Ucal}{\mathcal{U}}

\newcommand{\Bbf}{{\ensuremath{\mathbf{B}}}}

\newcommand{\Ibf}{{\ensuremath{\mathbf{I}}}}

\newcommand{\Mbf}{{\ensuremath{\mathbf{M}}}}

\newcommand{\Pbf}{{\ensuremath{\mathbf{P}}}}

\newcommand{\Sbf}{{\ensuremath{\mathbf{S}}}}
\newcommand{\Tbf}{{\ensuremath{\mathbf{T}}}}

\newcommand{\Vbf}{{\ensuremath{\mathbf{V}}}}

\newcommand{\Phbf}{{\ensuremath{\mathbf{\Phi}}}}

\newcommand{\zero}{{\ensuremath{\mathbf{0}}}}

\newcommand{\defing}[1]{\textbf{\emph{\mathversion{bold}#1}}}

\newcommand{\R}{\ensuremath{\mathbb{R}}}

\newcommand{\Cx}{\ensuremath{\mathbb{C}}}

\newcommand{\N}{\ensuremath{\mathbb{N}}}

\renewcommand{\leq}{\leqslant}

\newcounter{mcnt}

\newcounter{wordcnt}

%% file: mains/banach-povm.tex
\title[Banach--Ces\`{a}ro POVM]{Banach--Ces\`{a}ro positive operator-valued measure}

\author{Sviatoslav V. Dzhenzher}

\begin{abstract}
    We construct a family of bounded, finitely additive, covariant positive operator-valued measures (POVMs) on the real line by applying generalised Banach--Ces\`aro asymptotic averaging to locally finite covariant POVMs. In contrast to the strictly shift-invariant case, where countably additive operator measures are restricted to scalar multiples of a single fixed operator, the introduced covariant setting allows for non-commuting averaged operator values while successfully normalising the total mass. We develop an integration theory for complex-valued functions by taking the completion of simple integrable functions with respect to the continuous operator norm induced by the positive integral of their absolute values. The integration spaces for the original and the averaged measures are shown to possess distinct null classes. Finally, we introduce the operator-valued Fourier transform on these integration spaces and analyse the core algebraic constraints imposed by the covariance property on the exponential multipliers.
\end{abstract}

\thanks{\hspace{-5mm}
S.\,V. Dzhenzher: sdjenjer@yandex.ru. orcid: 0009-0008-3513-4312
\\
Moscow Institute of Physics and Technology 141701, Institutskii lane, 9, Dolgoprudny, Russia}

\maketitle
\thispagestyle{empty}


\vspace{3mm}

\section{Introduction}

Measures with values in Banach spaces are well studied \cite{DiestelUhl}.
They naturally arise when speaking of limit theorems for random operators \cite{Sakbaev-Shmidt-Shmidt-2022, DzhenzherSakbaev26-unbound}.
A physical motivation comes from the theory of open quantum systems \cite{Accardi-Volovich-Lu-2002-open-sys}, and in particular from random quantum channels \cite{Kholevo-qsys-chan-inf, Kholevo-qprobstat}.

Positive operator-valued measures (POVMs) provide a mathematical framework for describing quantum measurements beyond the idealised setting of projection-valued observables \cite{Nielsen-Chuang}. They connect the operator description of a quantum system with the probabilities of measurement outcomes and play an important role in the analysis of joint measurements.
Symmetries are particularly important in this framework: a transformation of the outcome space should be compatible with the corresponding unitary transformation of the system. This requirement leads to covariant POVMs, which express the symmetry of a measurement without requiring its operator values to remain unchanged under transformations of the outcomes \cite{Werner-2004, Amosov-covpovm}.

Although the standard definition of a POVM includes countable additivity, finitely additive structures also arise in extensions of quantum measurement theory. In particular, Srinivas \cite{Srinivas} considered finitely additive expectation-valued measures in connection with successive measurements of observables with continuous spectra; in that framework, the associated joint probabilities need not be countably additive.
A related motivation comes from averaging over non-compact symmetry groups.
Werner \cite{Werner-2004} used invariant means to reduce a joint-measurement uncertainty problem to the covariant case and explicitly addressed the possibility that the resulting operator-valued expectation functional carries weight at infinity.
These developments motivate the study of POVMs alongside the usual countably additive measures.

The relation between scalar measures and quantum states provides another motivation for studying asymptotic operator-valued measures.
Blecher and Weaver \cite{Blecher-Weaver} characterised the existence of singular countably additive states on \(\Bcal(\ell^2(\kappa))\) in terms of Ulam real-valued measurable cardinals. Their results connect the additivity properties of scalar probabilities with the normal and singular parts of quantum states on operator algebras (that is, positive functionals on $\Bcal(\Hcal)$, or, generally speaking, on $C^*$-algebras).
A dynamical aspect of this connection was studied in \cite{Dzhenzher26-qsing}, where averages of unitary shift conjugations with respect to scalar probability measures were used to construct singularising quantum channels. In the setting considered there, countable additivity of the averaging measure ensures preservation of countable additivity of states, while the additional assumption that the measure vanishes on singletons ensures the singularisation effect.

The present construction considers asymptotic averaging at the level of positive operator-valued measures.
Evaluation in a state then produces a scalar distribution, so one operator-valued measure describes the asymptotic distributions associated with different states simultaneously.
When the averaged measure $\mu$ has total mass \(\Ibf\), every normal quantum state \(\rho\) induces the probability charge \(A\mapsto\rho(\mu(A))\). This charge vanishes on bounded sets and is therefore purely finitely additive. The operator-valued construction also retains the covariance of the original measure and allows non-commuting values. These features connect asymptotic scalar probabilities with the operator structure of quantum observables.

In the scalar setting, Banach-type averaging provides a way to describe translation-invariant statistical behaviour on the whole real line without introducing a preferred location or a finite observation window.
Normalised averages over expanding intervals admit generalised limits even when ordinary averages fail to converge.
The resulting finitely additive measures retain information about asymptotic density while assigning zero mass to bounded sets.
Their integration theory and their relation to counting measure through the Fourier transform were studied by Dzhenzher and Sakbaev \cite{EDzhenzherSakbaev26}.
This raises the question of whether the same construction can be extended to operator-valued measures while retaining a genuinely non-commutative range.

The choice of symmetry is essential to this extension. As shown below, strict translation invariance of a locally finite, countably additive positive operator-valued measure forces its values to be scalar multiples of a single fixed positive operator.
Such a measure cannot have non-commuting values.
Covariance is therefore not merely an additional assumption: it is the appropriate replacement for scalar translation invariance when the measure takes operator values rather than scalar values. The resulting problem is to combine the asymptotic averaging underlying Banach measures with the operator transformation law of covariant observables.

The present paper develops this operator-valued Banach--Ces\`aro construction.
Starting from a locally finite covariant positive operator-valued measure on the real line, we obtain a bounded, finitely additive covariant measure by generalised averaging.
The construction does not require the original measure to have finite total mass, and an explicit example (Example~\ref{ex:non-com}) shows that the averaged measure can be normalised and have non-commuting values.
We also develop the integration of complex-valued functions with respect to the original and averaged measures and compare their classes of integrable functions. The aim is to connect the scalar theory of Banach measures with the operator-valued framework of covariant observables and to provide a basis for studying the corresponding Fourier transforms.

The paper is organised as follows.
Section~\ref{s:main} presents the main example (Example~\ref{ex:m-t}) of the operator-valued analogue of the Lebesgue measure, and constructs an operator-valued analogue of Banach-Ces\`{a}ro measure.
It is shown there that in order to obtain the non-trivial generalisation of the real-valued measures, one should replace the translation-invariance of measures with the covariance property; see Proposition~\ref{p:tr-inv} and Theorem~\ref{t:banach-covariant}.
Section~\ref{s:integral} develops integration with respect to these measures, with the main results being Theorem~\ref{t:norm} (for the general properties of POVM integration) and both Theorem~\ref{t:int-1} and Theorem~\ref{t:int-inf} for the constructed Banach--Ces\`aro POVM.
Finally, Section~\ref{s:concl} concludes the paper with directions for further research, including some analogues of the Fourier transform and the connection with quantum channels.

\section{Background and main results}\label{s:main}

In this section, let \(\Hcal := \Lcal_2(\R\to\Cx, \lambda)\) be a separable Hilbert space of functions \(u\colon\R\to\Cx\) such that \(\abs{u}^2\) is integrable w.r.t. the classical Lebesgue measure~\(\lambda\).
Denote by \(\Bcal(\Hcal)^+\) the positive cone of bounded operators on~$\Hcal$.
Denote by \(\{\Sbf_t\}_{t\in \R}\) the unitary group of shifts \(\Hcal\to\Hcal\) acting by
\[
    \Sbf_t u(y) = u(y+t).
\]
Let \(\Mcal\subset 2^\R\) be a $\sigma$-algebra of the sets that are measurable w.r.t. the classical Lebesgue measure~$\lambda$.
Let \(m\colon \Mcal\to\Bcal(\Hcal)^+\) be a (partially defined) locally-finite $\sigma$-additive absolutely continuous (w.r.t.~$\lambda$) POVM, such that for any \(t\in\R\),
\[
    m(A+t) = \Sbf_t^* m(A) \Sbf_t.
\]
The last property is called the \defing{covariance} of the POVM; for details, see, for example, \cite{Amosov-covpovm}.
By absolute continuity we mean that \(m(A)\xrightarrow[\lambda(A)\to 0]{}\zero\) for bounded \(A\in\Mcal\), which is, by the Pettis theorem \cite[Theorem~I.2.1]{DiestelUhl}, equivalent to the fact that \(m\) vanishes on sets of zero measure~$\lambda$.
By local finiteness we mean that \(m(A)\in\Bcal(\Hcal)^+\) for a bounded $A$, and by partial definiteness we mean that $m$ can be undefined on unbounded sets.
For convenience, we will sometimes write \(m(A)=+\infty\).
We presume that if \(m(A) \in \Bcal(\Hcal)^+\) and \(B\subset A\) is measurable, then \(m(B) \leq m(A)\).

Denote by \(\Ibf\in\Bcal(\Hcal)^+\) the \emph{identity operator}.

\begin{example}\label{ex:m-t}
    Fix a positive operator \(\Tbf\in\Bcal(\Hcal)^+\).
    Define
    \[
        m(A) = m_\Tbf(A) := \int_A \Sbf_x^*\Tbf\Sbf_x \,dx.
    \]
    Then $\sigma$-additivity and absolute continuity are clear.
    Local finiteness follows since
    \[
        \norm{m[\,a,b\,]} \leq \norm{\Tbf}\lambda[\,a,b\,].
    \]
    Moreover,
    \[
        \zero \leq m(A) \leq \norm{\Tbf}\lambda(A)\Ibf.
    \]
    Finally, the covariance holds since
    \[
        m(A+t) = \int_A \Sbf_{x+t}^*\Tbf\Sbf_{x+t} \,dx = \Sbf_t^* m(A) \Sbf_t.
    \]

    If we additionally demand that
    \[
        c\Ibf \leq \Tbf \leq C\Ibf
    \]
    for some \(C>c>0\), then
    \[
        c \lambda(A)\Ibf \leq m(A) \leq C\lambda(A)\Ibf.
    \]
    In particular, \(m(\R)\) is not a bounded positive operator.
\end{example}

\begin{proposition}\label{p:tr-inv}
    If the property of covariance in the definition of $m$ was changed to the shift invariance, then the measure $m$ would be of the form
    \[
        m(A) = \lambda(A) m(0,1).
    \]
\end{proposition}

\begin{proof}
    Indeed, for any \(u\in\Hcal\), we have that the measure
    \[
        \nu_u(A) :=(u, m(A) u)
    \]
    is non-negative, locally finite, absolutely continuous w.r.t.~$\lambda$, and shift invariant, and hence \(\nu_u(A) = \lambda(A)\nu_u(0,1)\).
\end{proof}

The following lemma shows that the covariance is a strong property, which implies absolute continuity.

\begin{lemma}
    Let \(m\colon \Mcal\to\Bcal(\Hcal)^+\) be a (partially defined) locally-finite $\sigma$-additive covariant POVM.
    Then $m$ is absolutely continuous w.r.t.~$\lambda$.
\end{lemma}

\begin{proof}
    Let \(N \in \Mcal\) be of \(\lambda(N)=0\).
    Due to local finiteness, $\sigma$-finiteness and $\sigma$-additivity, it is sufficient to consider $N$ being bounded.
    Then \(N + [\,-1,1\,] \subset [\,-x,x\,]\) for some \(x>1\).
    For \(u\in \Hcal\), consider the measure
    \[
        \nu_u(A) := (u, m(A)u).
    \]
    On the one hand, by the Fubini--Tonelli theorem,
    \[
        \int_{-1}^1 \nu_u(N+t)\,dt = \int_{-x}^xd\nu_u(y)\int_{-1}^1 \Ical_N(y-t)\,dt =
        \int_{-x}^x \lambda\bigl((y-N)\cap [\,-1,1\,]\bigr)\,d\nu_u(y) = 0,
    \]
    where the latter holds since \(\lambda(N)=0\) and since $\lambda$ is monotone and translation invariant.
    Thus, \(\nu_u(N+t)=0\) for $\lambda$-a.e. \(t\in[\,-1,1\,]\).
    On the other hand,
    \[
        \nu_u(N+t) = (\Sbf_t u, m(N)\Sbf_t u)
    \]
    is continuous by $t$, since \(\Sbf_t\) is strongly continuous on $\Hcal$.
    Hence, \(\nu_u(N+t) = 0\) for all \(t\in [\,-1,1\,]\); in particular, for $t=0$.
    Therefore, \(\nu_u(N) = 0\) for any $u\in \Hcal$, which provides \(m(N)=\zero\).
\end{proof}

It is probably worth noticing that the extension of a measure $m$ from an algebra of finite unions of segments can be done analogously to Proposition~\ref{p:tr-inv} for \(\nu_u(A) = (u,m(A)u)\), using the $\sigma$-finiteness and classical extension results for operator-valued measures \cite{DiestelUhl}.

Recall that a \emph{filter} is a proper subset of some powerset that is closed under intersections and taking supersets.
A \emph{ultrafilter} is a maximal filter.
A ultrafilter is said to be \emph{non-principal} if the intersection of all its elements is empty.
For details see, for example, \cite{Conway-lin-op, Dunford-Schwartz-vol1, EDzhenzherSakbaev26}.

Let \(\Ucal\subset 2^\R\) be a non-principal ultrafilter containing all neighbourhoods of \(+\infty\).
For \(A\in \Mcal\), define the \defing{Banach--Ces\`{a}ro POVM}
\[
    \mu_\Ucal(A) := \lim_{x\to\Ucal} \frac{m(A\cap [\,-x,x\,])}{2x}
\]
as the weak limit along the ultrafilter~$\Ucal$; in other words,
\[
    \left(u, \mu_\Ucal(A)v\right) = \lim_{x\to\Ucal} \frac{\left(u,m(A\cap [\,-x,x\,])v\right)}{2x};
\]
for the definition and properties of these limits, see, for example, \cite{EDzhenzherSakbaev26, Dunford-Schwartz-vol1, Conway-lin-op}.
Let \((\R, \Acal_\Ucal, \overline{\mu_\Ucal})\) be the completion of this measure.
In other words, the algebra \(\Acal_\Ucal\) consists of all sets of the form \(E \bigtriangleup N\), where \(E\in\Mcal\), and \(N\subset A\) for some \(\mu_\Ucal(A)=\zero\).
The completion is defined by
\[
    \overline{\mu_\Ucal}(E\bigtriangleup N) := \mu_\Ucal(E).
\]
Below, we will identify the original measure and the completion, so \(\mu_\Ucal\) is a finitely-additive complete measure.

\begin{proposition}\label{p:ultralim}
    The limit along the ultrafilter in the definition of \(\mu_\Ucal\) exists.
\end{proposition}

\begin{proof}
    Set \(C_0 := \norm{m[\,0,1\,]}\).
    Hence, from the covariance of $m$,
    \[
        \norm{m[\,-x,x\,]} \leq (2x+1)C_0,
    \]
    and \(\frac{m(A\cap [\,-x,x\,])}{2x}\) are bounded uniformly for $x>1$.
    It remains to use that any bounded function has the limit along the ultrafilter; see, for example, \cite{Dunford-Schwartz-vol1} for details.
\end{proof}

\begin{theorem}\label{t:banach-covariant}
    The measure \(\mu_\Ucal\) is covariant. That is, for any \(t\in\R\) and \(A\in\Acal_\Ucal\),
    \[
        \mu_\Ucal(A+t) = \Sbf_t^*\mu_\Ucal(A)\Sbf_t.
    \]
\end{theorem}

\begin{proof}
    Fix any $t\in\R$ and \(E\in\Mcal\).
    Then, for \(x \gg \abs{t}\),
    \[
        \norm{\frac{m((E+t) \cap [\,-x,x\,])}{2x} - \Sbf_t^*\frac{m(E\cap[\,-x,x\,])}{2x}\Sbf_t} \leq \frac{\norm{m[\,0,\abs{t}\,]}}{x} \xrightarrow[x\to +\infty]{} 0.
    \]
    Hence \(\mu_\Ucal(E+t) = \Sbf_t^*\mu_\Ucal(E)\Sbf_t\).

    Let \(N \subset A \in\Mcal\) be a subset of a measurable set of zero measure.
    Then \(N+t\subset A+t\), and
    \[
        \mu_\Ucal(N+t) \leq \mu_\Ucal(A+t) = \Sbf_t^*\mu_\Ucal(A)\Sbf_t = \zero.
    \]

    Now, for \(E \bigtriangleup N \in \Acal_\Ucal\) we have
    \[
        \mu_\Ucal(E \bigtriangleup N + t) = \mu_\Ucal((E+t) \bigtriangleup (N+t)) = \mu_\Ucal(E+t) = \Sbf_t^*\mu_\Ucal(E)\Sbf_t = \Sbf_t^*\mu_\Ucal(E \bigtriangleup N)\Sbf_t.\qedhere
    \]
\end{proof}

\begin{example}\label{ex:non-com}
    Here we show that, depending on \(m\), we are able to obtain some non-trivial non-commutative measures.

    Denote by \(\Fcal\colon\Hcal\to\Hcal\) the classical unitary Fourier transform.
    Let \(\Pbf_0\in\Bcal(\Hcal)\) be the orthogonal projector on the closed subspace of functions $f$ such that \(\supp \Fcal f \subset [\,0,1\,]\),
    and \(\Pbf_1\) be the analogous projector for \([\,1,2\,]\).
    Let \(\Mbf\in\Bcal(\Hcal)\) be the operator of multiplication by \(y\mapsto e^{iy}\).
    Set
    \[
        \Vbf := \Mbf\Pbf_0.
    \]
    Then
    \[
        \Sbf_t^*\Vbf\Sbf_t = e^{-it}\Vbf,
        \quad
        \Vbf^*\Vbf = \Pbf_0,
        \quad
        \Vbf\Vbf^* = \Pbf_1.
    \]
    Let
    \[
        \Tbf := \Ibf + \varepsilon (\Vbf+\Vbf^*)
    \]
    for \(0<\varepsilon<\frac{1}{2}\).
    Then
    \[
        m(A) = \lambda(A)\Ibf + \varepsilon\left(\Vbf \int_A e^{-ix}\,dx + \Vbf^*\int_A e^{ix}\,dx\right),
    \]
    and
    \[
        \mu_\Ucal(\R) = \Ibf.
    \]
    Moreover,
    \begin{gather*}
        \mu_\Ucal\{\cos > 0\} = \frac{1}{2}\Ibf + \frac{\varepsilon}{\pi}(\Vbf+\Vbf^*), \\
        \mu_\Ucal\{\sin > 0\} = \frac{1}{2}\Ibf - \frac{i\varepsilon}{\pi}(\Vbf-\Vbf^*).
    \end{gather*}
    Hence, the commutator is equal to
    \[
        [\mu_\Ucal\{\cos > 0\}, \mu_\Ucal\{\sin > 0\}] = \frac{2i\varepsilon^2}{\pi^2}(\Pbf_1 - \Pbf_0) \neq \zero.
    \]
\end{example}

\section{Integration}\label{s:integral}

For convenience, for a function $f\colon X\to Y$ and a set $B\subset Y$ we shorten
\[
    \{f\in B\}:= \{x\in X: f(x)\in B\} = f^{-1}(B).
\]

Here, let \(\Hcal\) be a Hilbert space.
Let \(\nu\colon\Mcal\to\Bcal(\Hcal)^+\) be a measure on an algebra or a $\sigma$-algebra $\Mcal$ of subsets of~$\Omega$.
Our examples of the main interest are the analogue $m$ of the Lebesgue measure, and the Banach--Ces\`{a}ro measure~\(\mu_\Ucal\).

We say that a function \(\phi\colon \Omega\to\Cx\) is \emph{simple} if it is $\nu$-measurable and has a finite number of values.
Denoting by
\[
    \Ical_A = \Ical(A)\colon y \mapsto
    \begin{cases}
        1, & y\in A, \\
        0, & y\notin A
    \end{cases}
\]
the \emph{indicator function} of a set $A\in\Mcal$,
for a simple function $\phi$ we now can write
\[
    \phi = \sum_{y \in \Cx} y \Ical\{\phi=y\}.
\]
Accordingly, for a locally-finite positive operator-valued measure $\nu$ we define
\[
    \int \phi\,d\nu := \sum_{y \in \Cx} y \cdot \nu\{\phi=y\} \in \Bcal(\Hcal),
\]
if all \(\nu\{\phi=y\}\) are defined; in other words, if \(\nu\{\phi\neq 0\}\) is defined.
As in the usual case of real-valued measures, in the sum above we assume \(0\cdot \infty = \infty\cdot \zero = \zero\).

If \(\int\phi\,d\nu\) exists, then \(\phi\) is said to be \defing{$\nu$-integrable}, or just integrable.
Denote by \(\Scal_1(\Omega\to\Cx,\nu)\) the linear vector space of simple integrable functions. Define the norm by
\[
    \norm{\phi}_\nu := \norm{\int \abs{\phi}\,d\nu}.
\]
Finally, we say that \(\phi\) and \(\psi\) are \defing{$\nu$-equivalent}, \(\phi\sim_\nu \psi\), if
\[
    \nu\{\phi\neq\psi\} = \zero.
\]

\begin{proposition}
    The function \(\phi \to \norm{\phi}_\nu\) is indeed the norm on the linear space \(\Scal_1(\Omega\to\Cx,\nu)\) (up to equivalence), and the integral is the positive linear operator.
\end{proposition}

\begin{proof}
    First, note that
    \[
        \{\phi+\psi \neq 0\} \subset \{\phi\neq 0\} \cup \{\psi\neq 0\},
    \]
    so if both \(\phi\) and \(\psi\) are integrable, then so is $\phi+\psi$.
    Hence, \(\Scal_1(\Omega\to\Cx,\nu)\) is indeed the linear space.

    The linearity of the integral follows by the finite additivity of $\nu$.
    The additivity follows from the linearity.
    The monotonicity follows by the positivity of $\nu$.

    The function \(\norm{\cdot}_\nu\) satisfies the triangle inequality since
    \[
        \norm{\phi+\psi}_\nu = \norm{\int\abs{\phi+\psi}\,d\nu} \leq 
        \norm{\int\abs{\phi}\,d\nu} + \norm{\int\abs{\psi}\,d\nu} =
        \norm{\phi}_\nu+\norm{\psi}_\nu.
    \]

    Finally, if \(\nu\{\phi\neq 0\} \neq 0\), then \(\norm{\phi}_\nu\neq 0\) by the positivity of~$\nu$.
\end{proof}

Now, let \(\Lcal_1(\Omega\to\Cx,\nu)\) denote the completion of the space of integrable simple functions.
The integral here is extended as a linear continuous operator, and the norm is defined by
\[
    \norm{f}_\nu = \norm{\int \abs{f}\,d\nu}.
\]
In order to relate the completion and the functions \(\Omega\to\Cx\), we need to say that a function \(f\colon\Omega\to\Cx\) is \defing{$\nu$-integrable} if there exists a Cauchy-\(\norm{\cdot}_\nu\) sequence \(\{\phi_n\}\) of simple integrable functions that converges to $f$ in $\nu$, \(\phi_n\xrightarrow[n\to\infty]{\nu} f\); the latter means that
\[
    \nu\bigl(\{\abs{f-\varphi_n}>\varepsilon\} \cap K\bigr) \xrightarrow[n\to\infty]{}\zero
\]
for any \(\varepsilon>0\) and any \(K\) of finite measure~$\nu$, where the convergence in the norm topology.

\begin{proposition}
    Let \(f\) be an integrable function \(\Omega\to\Cx\).
    Then
    \[
        \norm{\int f\,d\nu} \leq \norm{f}_\nu.
    \]
    If additionally $f$ is bounded and \(\nu\) is bounded, then
    \[
        \norm{f}_\nu \leq \norm{f}_\infty \norm{\nu(\R)}.
    \]
\end{proposition}

\begin{proof}
    For a finite function \(\phi = \sum_{y\in \Cx} y \Ical\{\phi=y\}\), for any \(u,v\in \Hcal\) we have
    \begin{multline*}
        \abs{\left(u, \left(\int \phi\,d\nu\right)v\right)} \leq \sum \abs{y}\abs{(u, \nu\{\phi=y\}v)} \leq \\
        \leq \sum \abs{y}\sqrt{(u, \nu\{\phi=y\}u)}\sqrt{(v, \nu\{\phi=y\}v)}
        \leq \sqrt{\left(u, \left(\int\abs{\phi}\,d\nu\right)u\right)}\sqrt{\left(v, \left(\int\abs{\phi}\,d\nu\right)v\right)},
    \end{multline*}
    where the last two inequalities are Cauchy--Bunyakovsky--Schwarz.
    Then, taking the supremum over \(\norm{u}=\norm{v}=1\),
    \[
        \norm{\phi\,d\nu} \leq \norm{\abs{\phi}\,d\nu}.
    \]
    After that, for an integrable function $f$, the inequality follows by the properties of completion.
    The additional part is clear for simple functions, and for arbitrary functions is proved by the properties of limits.
\end{proof}

\begin{lemma}\label{l:cauchy}
    Let \(\{\phi_n\}\colon \N\to\Scal_1(\Omega\to\Cx, \nu)\) be a $\norm{\cdot}_\nu$-Cauchy sequence of simple integrable functions, and \(\phi_n\xrightarrow[n\to\infty]{\nu}0\).
    Then \(\phi_n\xrightarrow[n\to\infty]{} 0\) in norm \(\norm{\cdot}_\nu\).
\end{lemma}

\begin{proof}
    Fix any \(\eta>0\).
    Since \(\phi_n\) is the Cauchy sequence, we can take \(N\) such that for any \(n,k>N\),
    \[
        \norm{\phi_n-\phi_k}_\nu < \eta.
    \]
    Set \(M_n := \max \Img\abs\phi_n\).
    Then for any \(\varepsilon>0\),
    \[
        \abs{\phi_n} \leq \abs{\phi_n-\phi_k} + \varepsilon\Ical\{\phi_n\neq 0\} + M_n\Ical(\{\phi_n\neq 0\}\cap \{\abs{\phi_k}>\varepsilon\}).
    \]
    Since all functions here are simple and integrable,
    \begin{multline*}
        \norm{\phi_n}_\nu \leq \norm{\phi_n-\phi_k}_\nu + \varepsilon \norm{\nu\{\phi_n\neq 0\}} + M_n\norm{\nu(\{\phi_n\neq 0\}\cap \{\abs{\phi_k}>\varepsilon\})} \leq \\ \leq
        \eta + \varepsilon \norm{\nu\{\phi_n\neq 0\}} + M_n\norm{\nu(\{\phi_n\neq 0\}\cap \{\abs{\phi_k}>\varepsilon\})}.
    \end{multline*}
    Now with $n,\varepsilon$ fixed, let \(k\to\infty\). Then from the convergence of $\phi_k$ to zero in measure,
    \[
        \norm{\phi_n}_\nu \leq \eta + \varepsilon \norm{\nu\{\phi_n\neq 0\}}.
    \]
    Since \(\varepsilon\) is arbitrary, we have
    \[
        \norm{\phi_n}_\nu \leq \eta.
    \]
    Since $\eta$ was chosen arbitrarily, this concludes the proof.
\end{proof}

We say that two functions \(f\) and \(g\) from $\Lcal_1(\Omega\to\Cx,\nu)$ are \defing{$\nu$-equivalent}, \(f\sim_\nu g\), if
for any \(\varepsilon>0\) and any \(K\) of finite measure $\nu$,
\[
    \nu(\{\abs{f-g}>\varepsilon\} \cap K) = \zero.
\]
For a $\sigma$-additive bounded measure $\nu$, it is equivalent to
\[
    \nu\{f\neq g\} = \zero.
\]

\begin{theorem}\label{t:norm}
    The function \(\norm{\cdot}_\nu\) on \(\Lcal_1(\Omega\to\Cx,\nu)\) is the norm (up to equivalence).
\end{theorem}

\begin{proof}
    The triangle inequality follows from the properties of the limit in the definition of the integral.

    Take any \(f\in\Lcal_1(\Omega\to\Cx, \nu)\).
    Since we are interested in \(\abs{f}\), w.l.o.g. we can assume that $f$ is non-negative.
    Take any \(\varepsilon>0\).
    Then for any \(K\) of finite measure $\nu$,
    \[
        \nu(\{f>\varepsilon\}\cap K) \leq \frac{1}{\varepsilon}\int f\,d\nu.
    \]
    Hence, if \(\norm{f}_\nu = 0\), then \(f\sim_\nu 0\).
    Vice versa, if \(f\sim_\nu 0\) and we additionally know that $\nu$ is bounded, then
    \[
        \int f\,d\nu \leq \varepsilon \nu(\R),
    \]
    and so \(\norm{f}_\nu = 0\).
    
    In case \(\nu(\Omega)\) is undefined, we need Lemma~\ref{l:cauchy}.
    Suppose that \(f\sim_\nu 0\).
    Take simple functions \(\{\phi_n\}\) that define the integrability of $f$.
    Since
    \[
        \{\abs{\phi_n} > 2\varepsilon\} \subset \{\abs{\phi_n-f}>\varepsilon\} \cup \{f > \varepsilon\},
    \]
    for any \(K\) of finite measure, we have
    \[
        \nu(K \cap \{\abs{\phi_n} > 2\varepsilon\}) \leq \nu(K \cap \{\abs{\phi_n-f} > \varepsilon\}) \xrightarrow[n\to\infty]{} \zero.
    \]
    Then \(\phi_n\xrightarrow[n\to\infty]{\nu} 0\), and by Lemma~\ref{l:cauchy}, \(\norm{\phi_n}_\nu\xrightarrow[n\to\infty]{}0\).
    Hence, \(\norm{f}_\nu = 0\).
\end{proof}

At last, we are ready to apply the obtained results for the Hilbert space \(\Hcal = \Lcal_2(\R\to\Cx, \lambda)\) and the measures \(m\) and \(\mu_\Ucal\).

\begin{theorem}\label{t:int-1}
    Let \(\mu_\Ucal\) be the Banach--Ces\`aro measure constructed by the measure $m$.
    Then
    \[
        \Lcal_1(\R\to\Cx,m) \subset \Lcal_1(\R\to\Cx,\mu_\Ucal).
    \]
    Moreover, every \(m\)-integrable function represents zero in \(\Lcal_1(\R\to\Cx,\mu_\Ucal)\).
\end{theorem}

\begin{proof}
    If \(\phi\in\Scal_1(\R\to\Cx,m)\), then \(\{\phi \neq 0\}\) is a set of finite measure~$m$, and so \(\mu_\Ucal\{\phi\neq 0\}=\zero\).
    So, for simple functions, the theorem is proved.
    
    Now, let \(f\in\Lcal_1(\R\to\Cx,m)\), and let \(\{\phi_n\}\) be a sequence of simple functions that represents the $m$-integrability of~$f$.
    Since $\{\phi_n\}$ is Cauchy-\(\norm{\cdot}_m\), we have
    \[
        C := \sup_{n} \norm{\phi_n}_m < +\infty.
    \]
    Hence, for any \(\varepsilon>0\) and \(x \gg 1\), from
    \[
        \{\abs{f} > 2\varepsilon\} \subset \{\abs{f-\phi_n} > \varepsilon\} \cup \{\abs{\phi_n} > \varepsilon\}
    \]
    we obtain
    \begin{multline*}
        m\bigl(\{\abs{f} > 2\varepsilon\} \cap [\,-x,x\,]\bigr) \leq
        m\bigl(\{\abs{f-\phi_n} > \varepsilon\} \cap [\,-x,x\,]\bigr) + \frac{1}{\varepsilon}\int\abs{\phi_n}\,dm \leq \\ \leq
        m\bigl(\{\abs{f-\phi_n} > \varepsilon\} \cap [\,-x,x\,]\bigr) + \frac{C}{\varepsilon}\Ibf.
    \end{multline*}
    With \(n\to\infty\) and \(\phi_n\xrightarrow[n\to\infty]{m} f\), we first summand in the latter vanishes, and
    \[
        \zero \leq m\bigl(\{\abs{f} > 2\varepsilon\} \cap [\,-x,x\,]\bigr) \leq \frac{C}{\varepsilon} \Ibf.
    \]
    So,
    \[
        \norm{\frac{m\bigl(\{\abs{f} > 2\varepsilon\} \cap [\,-x,x\,]\bigr)}{2x}} \leq \frac{C}{2x\varepsilon} \xrightarrow[x\to\Ucal]{} 0.
    \]
    Therefore,
    \[
        \mu_\Ucal\{\abs{f}>\varepsilon\} = 0,
    \]
    which means that \(\phi_n\), representing zeros in \(\Lcal_1(\R\to\Cx,\mu_\Ucal)\), are Cauchy-\(\norm{\cdot}_{\mu_\Ucal}\) and converge to $f$ in \(\mu_\Ucal\).
    Hence \(f \in \Lcal_1(\R\to\Cx,\mu_\Ucal)\) and \(\norm{f}_{\mu_\Ucal} = 0\).
\end{proof}

Denote by \(\Lcal_\infty(\R\to\Cx, m)\) the space of functions \(f\colon\R\to\Cx\) with the norm being the essential supremum.
The following is the generalisation of the analogous result from \cite{EDzhenzherSakbaev26}.

\begin{theorem}\label{t:int-inf}
    Let \(\mu_\Ucal\) be the Banach--Ces\`aro measure constructed by the measure $m$.
    Then
    \[
         \Lcal_\infty(\R\to\Cx,m) \subset \Lcal_1(\R\to\Cx,\mu_\Ucal).
    \]
    Moreover, for \(f\in\Lcal_\infty(\R\to\Cx,m)\),
    \[
        \int f\,d\mu_\Ucal = \lim_{x\to\Ucal} \frac{1}{2x} \int_{-x}^x f\,dm,
    \]
    where the limit is in the WOT sense.
\end{theorem}

\begin{proof}
    It is clear that a bounded function is \(\mu_\Ucal\)-integrable, since \(\mu_\Ucal(\R)\) is finite.
    It remains to notice that any function \(f\in\Lcal_\infty(\R\to\Cx, m)\) is bounded outside some set $A$ with \(m(A)=\zero\), and that \(\mu_\Ucal(A)=\zero\).

    Let us prove the `moreover' part.
    For an indicator, it is by definition; and for a simple function, it is by linearity.
    Now, let $f$ be an arbitrary function from \(\Lcal_\infty(\R\to\Cx,m)\).
    It is sufficient to consider the case when $f$ is bounded.
    So, since $f$ is bounded, there exist simple functions $\phi_n$ that converge to $f$ uniformly.
    As in the proof of Proposition~\ref{p:ultralim},
    there exists a constant $C_0>0$ such that
    \[
        \norm{\frac{1}{2x}\int_{-x}^x f\,dm - \frac{1}{2x}\int_{-x}^x \phi_n\,dm} \leq C_0 \norm{f-\phi_n}_\infty.
    \]
    Hence,
    \[
        \norm{\int f\,d\mu_\Ucal - \int \phi_n\,d\mu_\Ucal} \leq C_0 \norm{f-\phi_n}_\infty.
    \]
    Finally, for any \(u,v\in\Hcal\),
    \[
        \abs{\left(u, \left(\frac{1}{2x}\int_{-x}^x f\,dm - \int f\,d\mu_\Ucal\right)v\right)} \leq
        2C_0\norm{u}\norm{v}\norm{f-\phi_n}_\infty + \abs{
            \left(u, \left(\frac{1}{2x}\int_{-x}^x \phi_n\,dm - \int \phi_n\,d\mu_\Ucal\right)v\right)
        }.
    \]
    Taking \(n\to\infty\), we can make the first summand to be sufficiently small.
    After that it remains to take \(x\to\Ucal\).
\end{proof}

\section{Conclusions and further research}\label{s:concl}

We have constructed positive, finitely additive operator-valued measures by applying Banach--Ces\`aro averaging to locally finite covariant measures on the real line.
The construction preserves covariance with respect to the translation representation. In contrast to the strictly invariant, countably additive case, where the original measure has the form $\lambda(\cdot)\Tbf$, the covariant setting admits averaged measures with non-commuting values.
The example considered above shows that this non-commutativity is compatible with normalisation of the total mass.

The associated integration theory is obtained by completing simple functions in the norm induced by the positive integral of their absolute values. This provides a common framework for the locally finite original measure and the bounded, finitely additive averaged measure. Their integration spaces have different null classes: the identity on simple functions of finite $m$-mass induces the zero operator from $\Lcal_1(\R\to\Cx,m)$ to $\Lcal_1(\R\to\Cx,\mu_\Ucal)$, whereas $m$-essentially bounded functions define elements of $\Lcal_1(\R\to\Cx,\mu_\Ucal)$ through a bounded, generally
non-injective map.

A natural direction for further study is the Fourier transform
\[
    (\Fcal_\nu f)(y)
       =\int_\R f(x)e^{-ixy}\,d\nu(x),
       \qquad \nu\in\{m,\mu_\Ucal\}.
\]
The bounded-multiplier property of the constructed integration spaces makes
this definition meaningful for $f\in \Lcal_1(\R\to\Cx,\nu)$ and yields
\[
    \sup_{y\in\R}\norm{(\Fcal_\nu f)(y)} \leq \norm{f}_{1,\nu}.
\]
Thus $\Fcal_\nu$ takes complex-valued integrable functions to bounded $\Bcal(\Hcal)$-valued functions on the frequency axis.

For the averaged measure, set
\[
    C(r) := \int_\R e^{-irx}\,d\mu_\Ucal(x).
\]
Then the exponential $f_t(x):=e^{itx}$ satisfies
\[
    (\Fcal_{\mu_\Ucal}f_t)(y)=C(y-t),
\]
while covariance implies
\[
    \Sbf_a^*C(r)\Sbf_a=e^{ira}C(r).
\]
Unlike strict invariance, this relation does not force $C(r)$ to vanish for $r\ne0$. Consequently, the correspondence between exponentials and point masses established in the scalar setting (in \cite{EDzhenzherSakbaev26}) does not extend directly to arbitrary covariant operator-valued measures.

Further questions, therefore, concern the structure of $C$, the kernel and range of $\Fcal_\nu$, and the additional hypotheses under which a suitable Hilbert-space realisation admits an analogue of the scalar correspondence with counting measure.
The construction and integration results established here provide a starting point for these questions; no general inversion or unitary Fourier theorem is asserted.

Another possible direction of research is to try to replace the measure $m$ with some other measure, that, for example, is not covariant and is not absolutely continuous w.r.t.~$\lambda$.
For example, instead of \(\Hcal = \Lcal_2(\Rcal\to\Cx, \lambda)\), one could take the Hilbert space of square-integrable functions w.r.t. the counting measure \(A\mapsto \abs{A}\).
It might be interesting to find out the criterion for when the resulting Banach--Ces\`aro POVM remains covariant, even if $m$ is not.
See \cite{Beketov-26} for analogous examples of the proposed questions in the scalar setting.

Finally, the measures in Example~\ref{ex:m-t} also provide a direct connection with quantum channels in the Heisenberg picture. For \(\Bbf\in \Bcal(\Hcal)\), define
\[
    \Phbf_\Ucal(\Bbf) := \lim_{x\to\Ucal} \frac{1}{2x}\int_{-x}^{x}\Sbf_t^*\Bbf \Sbf_t\,dt.
\]
By the results from \cite{EDzhenzherSakbaev26}, it is clear that \(\Phbf_\Ucal\) is the WOT integral over the scalar Banach--Ces\`aro measure.
This \(\Phbf_\Ucal\) is a unital completely positive map, since these properties hold for the finite averages and are preserved by the bounded weak operator limit. For the measure generated by a positive operator \(\Tbf\), its averaged total mass satisfies
\[
    \mu_\Ucal^\Tbf(\R)=\Phbf_\Ucal(\Tbf).
\]
Following the ideas (and definitions) from \cite{Dzhenzher26-qsing}, it is natural to find out whether the quantum states \(\rho \circ \Phbf_\Ucal\) are singular for such measures for any quantum state \(\rho\in\Bcal(\Hcal)^*\).
This would provide a more physical construction of singularising quantum channels than what was proposed in \cite{Dzhenzher26-qsing}.

\section*{Declarations}

\subsection*{Conflict of interest}

There is no conflict of interest.

\subsection*{Data availability}

Not applicable.

\subsection*{Funding}

No funding was received during the research.

%% file: bibtex.bib
@article{Dzhenzher26-qsing,
    author = {S.~V.~Dzhenzher},
    title = {Singularising preserving countable additivity quantum channels on quantum measurable cardinals},
    journal = {Ann. Funct. Anal.},
    volume = {18},
    number = {1},
    pages = {4},
    year = {2027},
    doi = {10.1007/s43034-026-00551-8}
}

@article{DzhenzherSakbaev26-unbound,
    author = {S.~V.~Dzhenzher and V.~Zh.~Sakbaev},
    title = {The Strong Law of Large Numbers for random semigroups with unbounded generators on uniformly smooth Banach spaces},
    journal = {Infin. Dimens. Anal. Quantum Probab. Relat. Top.},
    year = {2026},
    %volume = {29},
    doi = {10.1142/S0219025726500049}
}

@article{Amosov-covpovm,
    author = {Amosov, G.},
    title = {On Constructing Informationally Complete Covariant Positive Operator-Valued Measures},
    journal = {Entropy (Basel)},
    year = {2023},
    volume = {25},
    issue = {5},
    doi = {10.3390/e25050783}
}

@article{Beketov-26,
    author = {Dmitriy I. Beketov},
    title = {Generalized Limits of Finitely Additive Measures on Locally Compact Groups},
    usera = {2611.21199},
    year = {2026}
}

@article{Blecher-Weaver,
    title = {Quantum measurable cardinals},
    journal = {J. Funct. Anal.},
    volume = {273},
    number = {5},
    pages = {1870-1890},
    year = {2017},
    issn = {0022-1236},
    doi = {10.1016/j.jfa.2017.05.006},
    author = {David P. Blecher and Nik Weaver}
}

@article{EDzhenzherSakbaev26,
    author = {E.~A.~Dzhenzher and V.~Zh.~Sakbaev},
    title = {Banach limits, Banach measure, and properties of spaces of integrable functions (in Russian)},
    journal = {Vladikavkaz. Math. J.},
    volume={28},
    issue={2(110)},
    pages={50--69},
    year = {2026},
}

@article{Sakbaev-Shmidt-Shmidt-2022,
author={V.~Zh.~Sakbaev.
and E.~V.~Shmidt
and V.~Shmidt},
title={Limit Distribution for Compositions of Random Operators},
journal={Lobachevskii J. Math.},
year={2022},
month=jul,
day={01},
volume={43},
number={7},
pages={1740-1754},
issn={1818-9962},
doi={10.1134/S199508022210033X}
}

@article{Srinivas,
    author = {M. D. Srinivas},
    title = {Collapse postulate for observables with continuous spectra},
    journal = {Commun. Math. Phys.},
    year = {1980},
    volume = {71},
    pages = {131--158},
    doi = {10.1007/BF01197917}
}

@article{Werner-2004,
    author = {R. F. Werner},
    title = {The uncertainty relation for joint measurement of position and momentum},
    journal = {Quantum Inf. Comput.},
    vol = {4},
    pages = {546--562},
    year = {2004}
}

@Inbook{Accardi-Volovich-Lu-2002-open-sys,
author="Accardi, Luigi
and Volovich, Igor
and Lu, Yun Gang",
title="Open Systems",
bookTitle="Quantum Theory and Its Stochastic Limit",
year="2002",
publisher="Springer Berlin Heidelberg",
address="Berlin, Heidelberg",
pages="113--151",
isbn="978-3-662-04929-7",
doi="10.1007/978-3-662-04929-7_4"
}

@Inbook{Conway-lin-op,
author="Conway, John B.",
title="Linear Operators on a Banach Space",
bookTitle="A Course in Functional Analysis",
year="2007",
publisher="Springer New York",
address="New York, NY",
pages="166--186",
isbn="978-1-4757-4383-8",
doi="10.1007/978-1-4757-4383-8_6"
}

@book{DiestelUhl,
  title={Vector Measures},
  author={Diestel, J. and Uhl, J.~J.},
  isbn={9780821815151},
  lccn={lc77009625},
  series={Math. Surv. Monogr.},
  url={https://books.google.ru/books?id=fEKZAwAAQBAJ},
  year={1977},
  publisher={AMS}
}

@book{Dunford-Schwartz-vol1,
  title={Linear Operators. Part 1. General Theory},
  author={Dunford, N. and Schwartz, J.T.},
  series={Pure and Applied Mathematics},
  url={https://books.google.ru/books?id=x94XzQEACAAJ},
  year={1958},
  publisher={Interscience Publ.},
  address={New York}
}

@book{Kholevo-qprobstat,
author       ={A.~S. Kholevo },  
title        ={Quantum probability and quantum statistics },  
publisher    ={Probability theory -8, Itogi Nauki i Tekhniki. Series Sovrem. Probl. Mat. Fund. Napr.},    
volume       ={83} ,    
address      ={VINITI, Moscow },    
pages        ={5--132 },    
year         ={1991 }    
}

@book{Kholevo-qsys-chan-inf,
  title={Quantum Systems, Channels, Information: A Mathematical Introduction},
  author={Kholevo, A.~S.},
  isbn={9783110273250},
  lccn={2012033596},
  series={De Gruyter studies in mathematical physics},
  doi={10.1515/9783110273403},
  %url={https://books.google.ru/books?id=sObHMgEACAAJ},
  year={2012},
  publisher={De Gruyter}
}

@book{Nielsen-Chuang,
    place={Cambridge},
    title={Quantum Computation and Quantum Information: 10th Anniversary Edition},
    publisher={Cambridge University Press},
    author={Nielsen, Michael A. and Chuang, Isaac L.},
    year={2010}
}
